\documentclass[12pt]{article}       
\usepackage{array}
\usepackage{amsmath}
\usepackage{amssymb}
\begin{document}

\newcommand{\comment}[1]{}    
\newcommand{\hs}{\enspace}
\newcommand{\hhs}{\thinspace}
\newcommand{\real}{\ifmmode {\rm R} \else ${\rm R}$ \fi}
\def\nat{\mbox{\vrule height 7pt width .7pt depth 0pt\hskip -.5pt\bf N}}
\newcommand{\qed}{\hfill{\setlength{\fboxsep}{0pt}
                  \framebox[7pt]{\rule{0pt}{7pt}}} \newline}

\newcommand{\eqed}{\qquad{\setlength{\fboxsep}{0pt}
                  \framebox[7pt]{\rule{0pt}{7pt}}}\newline }
\newtheorem{theorem}{Theorem}
\newtheorem{lemma}[theorem]{Lemma}         
\newtheorem{corollary}[theorem]{Corollary}
\newtheorem{definition}[theorem]{Definition}
\newtheorem{claim}[theorem]{Claim}%
\newtheorem{conjecture}[theorem]{Conjecture}
\newtheorem{proposition}[theorem]{Proposition}
\newtheorem{construction}[theorem]{Construction}
\newtheorem{problem}[theorem]{Problem}

\newcommand{\proof }{{\bf Proof: }}          



\newcommand{\seq}{d_1, \ldots, d_{2k+2}}
\def\eps{\varepsilon}

\def\n{\notag}

\def\OO{\mathcal{O}}
\def\P{P\mathcal{P}}
\def\GG{\mathcal{G}}
\def\AA{\mathcal{A}}
\def\BB{\mathcal{B}}
\def\HH{\mathcal{H}}
\def\CC{\mathcal{C}}
\def\QQ{\mathcal{Q}}
\def\e{\varepsilon}
\def\FF{\mathcal{F}}


\title{Specified Intersections}
\author{
 Dhruv Mubayi
\thanks{Department of Mathematics, Statistics, and Computer
Science, University of Illinois,  Chicago, IL 60607. Research
supported in part by NSF grant DMS-0969092. Email: {\tt
mubayi@math.uic.edu}} \quad and \quad Vojtech R\"odl \thanks{Department
of Mathematics and Computer Science, Emory University.  Research supported in part by
NSF grant DMS-0800070. Email: {\tt rodl@mathcs.emory.edu}} }

\date{\today}
\maketitle

\vspace{-0.3in}

\begin{abstract}
Let $M \subset [n]:=\{0, \ldots, n\}$ and ${\AA}$ be a family of subsets of an $n$ element set such that $|A \cap B| \in M$
for every $A, B \in \AA$.
Suppose that $l$ is the maximum number  of consecutive integers contained in $M$  and $n$ is sufficiently large. Then
$$|\AA|< \min\{1.622^n 10^{2l+5} \, ,\quad 2^{n/2+l\log^2 n}\}.$$
 The first bound complements the previous bound of
roughly $(1.99)^n$ due to Frankl and the second author~\cite{FR} proved under the assumption that $M=[n]\setminus \{n/4\}$. For  $l=o(n/\log^2 n)$, the second bound above becomes better than the first bound.  In this case, it yields $2^{n/2+o(n)}$ and this
can be viewed as a generalization (in an asymptotic sense) of the famous Eventown theorem of Berlekamp~\cite{B} and Graver~\cite{G}. We conjecture that our bound $2^{n/2+o(n)}$ remains valid as long as $l<n/10$.

 Our second result complements the result of \cite{FR} in a different direction.  Fix $\eps>0$ and $\eps n < t<n/5$ and let $M=[n]\setminus (t, t+n^{0.525})$.  Then, in the notation above, we prove that for $n$ sufficiently large,

 $$|\AA| \le n{n \choose (n+t)/2}.$$
 This is essentially  sharp aside from the multiplicative factor of $n$. The short proof uses the Frankl-Wilson theorem and results about the distribution of prime numbers.
  We conjecture that a similar bound holds for $M=[n]\setminus \{t\}$ whenever $\eps n < t < n/3$.  A similar conjecture when $t$ is fixed and $n$ is large was earlier made by Frankl~\cite{F} and proved by
  Frankl and F\"uredi~\cite{FF}.

\end{abstract}

\section{Introduction}

Throughout this paper, we let $[n]:=\{0, 1, \ldots, n\}$ and $V$ denote an $n$-element set. Say that a family of sets $\AA$ is $M$-intersecting if for every $A, B \in \AA$, we have $|A \cap B| \in M$.  Suppose that $\AA \subset 2^{V}$. Our starting point is the following result.

\begin{theorem} \label{fr} {\bf (Frankl-R\"odl \cite{FR})}
For every $0<\eta<1/4$ there exists $\e>0$ and $n_0$ such that
if $n>n_0$, $\eta n <t < (1/2-\eta)n$, $M=[n]\setminus \{t\}$ and $\AA$
is $M$-intersecting, then
$|\AA|<(2-\e)^n$.
\end{theorem}
 Theorem \ref{fr} was previously
conjectured by Erd\H os, and  has applications in geometry~\cite{LR}, combinatorics~\cite{ESz}, coding theory, communication complexity~\cite{Sg} and quantum computing~\cite{BCW}. In words, the result says that if we forbid even one number $t\in (cn, (1/2-c)n)$ which is constant times $n$ away from both 0 and $n/2$
as an intersection size, then the size of our family must be exponentially smaller than the family of all sets. The result of \cite{FR} was actually more general.
Say that a pair $(\AA, \BB)$ of set systems is $M$-{\it intersecting}
if for every $A \in \AA$ and $B \in \BB$, we have $|A \cap B| \in M$.  Frankl and R\"odl proved that in the setup above, we have $|\AA||\BB|<(4-\e)^n$.  This is stronger, since we may let $\AA=\BB$.

At the other end of the spectrum, \cite{FR} also  proves that if $t \in [n]$ and
$M=\{t\}$, then
\begin{equation} \label{easyfr} |\AA||\BB|\le 2^n\end{equation}
 and this is sharp for many values of $t$.

In this paper we consider the size of $M$-intersecting families for two different types of $M$ which are in between these two extremes.

\subsection{Forbidding syndetic sets}
A set of integers is called $l$-syndetic if it intersects every interval of length $l$.  Also, for a set $M$ of integers, we define the length $l(M)$ to be the maximum number of consecutive integers contained in $M$.  Clearly $l(M)\le l$ iff $\overline M$ is $(l+1)$-syndetic.

 Our first result is concerned with finding upper bounds for families $\AA \subset 2^V$ that are
$M$-intersecting in terms of $l(M)$.  As $l(M)$ gets smaller (i.e. the forbidden set of intersection sizes intersects every interval of smaller length), this places more restrictions on $\AA$ and we therefore expect a better upper bound. Hence it is not surprising that as $l(M)$ becomes smaller, our bound is numerically better than the bound obtained in \cite{FR} for $M=[n]\setminus \{t\}$. As in \cite{FR}, we prove our result for pairs of families.

\begin{theorem} \label{main}
Let $M\subset [n]$ with $l(M)=l$.
Suppose that $(\AA, \BB)$ is an $M$-intersecting pair of families in $2^V$.  Then
$$|\AA||\BB| < \min\left\{2.631^n \times 10^{4l+10}, \quad 2^{n+2l\log^2n}\right\}.$$
\end{theorem}

{\bf Remarks.}

1) The constant $10^{4l}$ above has not been optimized and can be improved to slightly less than $10^{3l}$.

2)  The theorem is meaningful only for small $l$, say $l < n/10$.  Indeed, one quickly notices that if $l$ is a bit larger, say $l=0.15n$, then both bounds in the minimization are larger than $4^n$ (for large $n$) which is a trivial bound.  Therefore,  when $l>0.15n$, Theorem \ref{main} says nothing nontrivial.
 For this case upper bounds of the form $|\AA||\BB|<(4-\eps)^n$ follow only from Theorem \ref{fr} and a result of Sgall~\cite{Sg}.  When the two intersection sizes $n/3$ and $n/5$ are forbidden, the best upper bound is due to Sgall~\cite{Sg}.

\begin{corollary} \label{cor}
Let $M\subset [n]$ with $l(M)=l$.
Suppose that $\AA\subset 2^V$ is an $M$-intersecting family.  Then
$$|\AA|<\min\left\{1.622^n \times 10^{2l+5} \, , \quad 2^{n/2+l\log^2 n}\right\}.$$
\end{corollary}

{\bf Remarks.}

1) If $l \gg n/\log^2 n$, then the first bound in Corollary \ref{cor} is better and if $l \ll n/\log^2 n$, then the second bound in Corollary \ref{cor} is better.

2) The first bound in Theorem \ref{main} and Corollary \ref{cor} applies even when $l$ is linear in $n$, for example, when $l=n/10^4$ we get the upper bound $1.63^n$ from Corollary \ref{cor}. In this case, the forbidden set of intersection sizes $P=[n]\setminus M$ could have only $10^4$ numbers  in $[n]$ that are close to being uniformly distributed. We are not aware of any result that addresses such cases directly. As mentioned above, the only nontrivial bound we know follows directly from Theorem \ref{fr} and is about $(1.99)^n$ (though by carefully going through the calculations from \cite{FR} one could perhaps improve this slightly).  A related result of Sgall \cite{Sg} gives better bounds than Theorem \ref{fr} in the case when more than one intersection size is omitted, though the omitted sizes must correspond to congruence classes modulo some integer; Theorem \ref{main} does not require this.

3) If $l=o(n/\log^2n)$, then the second bound in Corollary \ref{cor} is $2^{n/2+o(n)}$ and  this can be viewed as a generalization (in an asymptotic sense) of the famous Eventown theorem of Berlekamp~\cite{B} and Graver~\cite{G}, which states that if $M=\{0,2,..\}$ and $|A|$ is even for every $A \in \AA$, then $|A| \le 2^{n/2}$. In particular, our bound does not require $M$ to be any collection of residue classes. Moreover, the second bound in Theorem \ref{main} can also be viewed as a generalization (when viewed as an asymptotic result) of (\ref{easyfr}) which applies for $M=\{t\}$ since we may let $l=1$.

The second bound in Theorem \ref{main} cannot be extended to $2^{n+o(n)}$ (independent of $l$) even for small values of $l$. For example,  one can let $l=n/10^4$, $M=[n]\setminus \{n/10^4, 2n/10^4,
\ldots, n\}$, $\AA=2^V$ and $\BB={[n] \choose \le n/10^4-1}$.  Then
$|A \cap B|\in M$ for every $A \in \AA$ and $B \in \BB$ and yet
$|\AA||\BB|>2^{1.0001n}$.
On the other hand, we are not able to obtain a construction of this type for just one family, and the best construction we have for Corollary \ref{cor} is obtained by the Eventown construction: assuming $n$ is even, take all subsets of $[n/2]$ and then double each point.  The resulting family has size $2^{n/2}$ and every two sets have even intersection size.  This leads us to make the following conjecture.

\begin{conjecture} \label{conj}
Let $1<l<n/10$ and $M\subset [n]$ with $l(M)=l$.
Suppose that $\AA\subset 2^V$ is an $M$-intersecting family.  Then
$$|\AA|<2^{n/2+o(n)}.$$
\end{conjecture}

{\bf Remark.}  The condition $l<n/10$ in Conjecture \ref{conj} is somewhat arbitrary, though some bound on $l$ of this type is required to prohibit constructions of the form ${[n] \choose \le l-1}$ with $M=[n]\setminus \{l, 2l, ..\}$.
Such constructions have larger size than $2^{n/2}$ if $l$ is large.
\bigskip

\subsection{Small intervals}
Our second result considers the case when $M$ omits a very small interval. In this case we prove an essentially sharp result for the maximum size of an $M$-intersecting family.
The starting point of this line of research is perhaps Katona's theorem~\cite{Kat} which determines the maximum size of an $M$-intersecting family of subsets of $[n]$ when $M=[n]\setminus [t]$; in other words, every two sets have at least $t+1$ elements in common. To state Katona's result precisely, define
$\AA(n,t)$ to be $\{A \subset V:|A| \ge (n+t+1)/2\}$ if $n+t$ is odd
and $\{A \subset V:|A \cap (V\setminus \{v\})| \ge (n+t)/2\}, v \in V$ is fixed,  if $n+t$ is even.

\begin{theorem} {\bf (Katona \cite{Kat})}
Let $\AA \subset 2^V$ and suppose that $|A \cap A'|>t$ for every $A, A' \in \AA$. Then
$$|\AA| \le |\AA(n,t)|.$$
Moreover, if $t \ge 1$ and $|\AA|=|\AA(n,t)|$, then $\AA=\AA(n,t)$.
\end{theorem}

The bound in Katona's theorem is essentially ${n \choose (n+t+1)/2}$, achieved by taking all large enough sets. If we weaken the hypothesis in Katona's theorem by forbidding just one intersection size, namely $t$, then Erd\H os~\cite{E} asked how large $|\AA|$ could be.
Later Frankl~\cite{F} conjectured that for $n>n_0(t)$,
$$|\AA| \le |\AA^*(n,t)|$$
 where $\AA^*(n,t)$ is obtained from $\AA(n,t)$ by adding all sets of size less than $t$. This was later proved by Frankl and F\"uredi for fixed $t$ (see also \cite{FFforb} for related results).

 \begin{theorem} \label{ff} {\bf (Frankl-F\"uredi \cite{FF})}
 Let $\AA \subset 2^V$ and suppose that $|A \cap A'|\ne t$ for every $A, A' \in \AA$. Then for $n>n_0(t)$,
$$|\AA| \le |\AA^*(n,t)|$$
and equality holds only for $\AA=\AA^*(n,t)$.
 \end{theorem}

  The condition $n>n_0(t)$ above appears to be essential in the argument of \cite{FF} and if we do not assume this, then the bound obtained from the proof in \cite{FF} is larger than $|\AA^*(n,t)|$ when $t$ is linear in $n$.
  

  In our final result we weaken the condition $n>n_0(t)$ to $n>5t$, but  enlarge the set of missing intersection sizes from one number (namely $t$) to a small interval  around $t$.  Under these conditions, we obtain an upper bound that is not exactly $|\AA^*(n,t)|$, though the logarithm of our upper bound is asymptotically equal to
  $\log |\AA^*(n,t)|$.

\begin{theorem}  \label{misint}
 Let $0<\eps<1/5$ be fixed, $n>n_0(\eps)$, $\eps n <t <n/5$ and $M=[n]\setminus (t, t+n^{0.525})$.
Suppose that $\AA$ is an $M$-intersecting family of subsets of $[n]$. Then
$$|{\AA}|< n{n \choose (n+t)/2}.$$
\end{theorem}

{\bf Remark.}  The constant 0.525 that appears above is a direct consequence of the result of Baker-Harman-Pintz \cite{BHP} that there is a prime in every interval $(s-s^{0.525}, s)$ as long as $s$ is sufficiently large.  


  We conjecture that $|\AA| < {n \choose (n+t)/2}2^{o(n)}$ for all $t <n/3$ even in the case when $M=[n]\setminus \{t\}$.  For $n/3\le t < (1/2-\e)n$, we conjecture that $|\AA| < {n \choose t} 2^{o(n)}$.

\section{Proof of Theorem \ref{main}}

We prove Theorem \ref{main} in two sections, each devoted to one of the bounds in the minimum.

\subsection{The first bound}
In this section we prove that $|\AA||\BB|<2.631^n 10^{4l+10}$.

{\bf Definitions and Notation.}

It is more convenient to phrase our proof in terms of complements of $M$, so we say that
$(\AA, \BB)$ is $P$-omitting if $|A \cap B| \not\in P$ for each $A \in \AA, B \in \BB$.
We will assume that $P=[n]\setminus M$ and $l(M)=l$ in the rest of this section.

Let $\FF \subset 2^V$ such that $\FF$ is $P$-omitting. Define
$$p(\FF)=\frac{|\FF|}{2^{|V|}}.$$
For all $v \in V$ let
$$\FF_1(v)=\{F\setminus \{v\}: v \in F \in \FF\} \quad \hbox{ and }
\quad \FF_0(v)=\{F \in \FF: v \not\in F\}.$$
 Note that $\FF_1(v), \FF_0(v) \subset  2^{V\setminus \{v\}}$.

Given a set $S\subset [n]$ and an integer $r$ let $S-r=\{s-r:s \in S\}$.

We now begin the proof of the first bound in Theorem~\ref{main}.
Let $\eps=2/10^4$ and $a_0=6-2\sqrt{5}<2$.  Put
$$f(x)=2-\frac{x}{2-\sqrt x}.$$
An easy calculation shows that $f(a_0)=0$.  Moreover, if $x$ is slightly less than $a_0$, then $f(x)>0$.  Hence we may choose
 $\delta=0.007$ such that
$f(a)>\sqrt\eps$, where $a=a_0-\delta$.

Before embarking on the proof,
let us make some preliminary observations.  Suppose that $\FF, \GG \in 2^V$ and $(\FF, \GG)$ is  $P$-omitting. Let $v \in V$ and write $\HH_i$ for $\HH_i(v)$ where $\HH \in \{\FF, \GG\}$ and $i \in \{0,1\}$.
Then

$\bullet$ $(\FF_1, \GG_1)$ is $P'$-omitting where $P'=P-1$

$\bullet$ $(\FF_0, \GG_0 \cup \GG_1)$ is $P$-omitting

$\bullet$ $(\FF_1, \GG_0 \cap \GG_1)$ is $P'$-omitting, where $P'=(P-1)\cup P$.

The most salient of the three properties above is the last one, since it implies that if $l([n]\setminus P)=l$, then $l([n]\setminus P')=l-1$.

The proof of the result, which extends the approach taken in \cite{FR}, is algorithmic.  Given a pair $(\FF, \GG)$ that is $P$-omitting where $l([n]\setminus P)=l$, we decompose it into the three pairs above. We will argue that the product of at least one of them must be large if $|\FF||\GG|$ is large.  In the first two cases, the families become more dense, while  the third case  when the family gets a bit sparser may happen only a few times.
\bigskip

{\bf Procedure.}

Recall that $\eps=2/10^4$, $a_0=6-2\sqrt{5}$, $\delta=0.007$ and $a=a_0-\delta$. 

{\bf Input:} A 4-tuple $(\FF, \GG, P, V)$ such that the
pair $(\FF, \GG)$  is $P$-omitting, $\FF, \GG \subset 2^V$ and $\FF, \GG \ne \emptyset$.

Suppose that there exists $v \in V$ such that for $\FF_{i}=\FF_{i}(v), \GG_{i}=\GG_{i}(v)$, $i=0,1$, one of the three possibilities (\ref{star}), (\ref{2star}), (\ref{3star}) below holds.

1) If
\begin{equation} \label{star}
p(\FF_1)p(\GG_1)>a\, p(\FF)p(\GG)
\end{equation}
holds, then  set $\FF'=\FF_1$ and $\GG'=\GG_1$ and repeat the procedure with $(\FF, \GG, P, V)$ replaced by $(\FF', \GG', P-1, V \setminus \{v\})$.

2) If (\ref{star}) fails but
\begin{equation} \label{2star}
p(\GG_0 \cup \GG_1)p(\FF_0)> a\, p(\FF)p(\GG)\end{equation}
holds,  then set $\FF'=\FF_0$ and $\GG'=\GG_0 \cup \GG_1$ and repeat the procedure with $(\FF, \GG, P, V)$ replaced by $(\FF', \GG', P, V \setminus \{v\})$.

3) If (\ref{star}) and (\ref{2star}) fail but
\begin{equation} \label{3star}
p(\GG_0 \cap \GG_1) p(\FF_1) >
\eps \, p(\FF) p(\GG) \end{equation}
holds, then set  $\FF'=\FF_1$ and $\GG'=\GG_0 \cap \GG_1$, $P'=(P-1) \cup P$ and repeat the procedure with $(\FF, \GG, P, V)$ replaced by $(\FF', \GG', P', V \setminus \{v\})$.

If (\ref{star}), (\ref{2star}) and (\ref{3star}) all fail, then {\bf Stop}.

Suppose that for all $v$, all of (\ref{star}), (\ref{2star}) and (\ref{3star}) fail and we have stopped the algorithm. Then
$$p(\FF_1)p(\GG_1)\le a\, p(\FF)p(\GG)$$
and hence for each $v \in V$, either
\begin{equation} \label{bar}p(\FF_1) \le \sqrt a \,p(\FF) \qquad \hbox{ or } \qquad
p(\GG_1) \le \sqrt a \,p(\GG) \end{equation}
Moreover, since (\ref{2star}) and (\ref{3star}) fail, we may assume
\begin{equation} \label{2bar}
p(\GG_0 \cup \GG_1)p(\FF_0)\le  a\, p(\FF)p(\GG)\end{equation}
and
\begin{equation} \label{3bar}
p(\GG_0 \cap \GG_1) p(\FF_1) \le
\eps \, p(\FF) p(\GG). \end{equation}
We will show that under the assumptions (\ref{bar}), (\ref{2bar}) and (\ref{3bar}), the inequalities in (\ref{bar}) can be strengthened as follows:

{\bf Claim.}  For all $v \in V$, either

a)  $p(\FF_1) \le \sqrt{\eps} p(\FF)$, or

b) $p(\GG_1) \le \sqrt{\eps} p(\GG)$.

{\bf Proof of Claim.}   We will use the identities
$$p(\HH_1)+p(\HH_0)=2p(\HH)=p(\HH_0 \cup \HH_1) + p(\HH_0 \cap \HH_1)$$ for $\HH \in \{\FF, \GG\}$. Let $v \in V$ be a vertex for which the first inequality in (\ref{bar}) holds.  Then
\begin{align}
\eps \,  p(\GG) p(\FF)\ge p(\GG_0 \cap \GG_1) p(\FF_1)
&=(2p(\GG) -p(\GG_0 \cup \GG_1))p(\FF_1) \n\\
&\ge \left(2p(\GG)-\frac{a\, p(\GG)p(\FF)}{p(\FF_0)}\right)p(\FF_1) \n \\
&=
 \left(2-\frac{a \, p(\FF)}{2p(\FF)-p(\FF_1)}\right)p(\GG)p(\FF_1) \n \\
&\ge \left(2-\frac{a}{2-\sqrt a}\right)p(\GG)p(\FF_1). \n \\
&\ge \sqrt{\eps} \, p(\GG) p(\FF_1). \n
\end{align}
Comparing the LHS and RHS of the inequalities above yields $p(\FF_1)\le \sqrt{\eps} p(\FF)$.  The proof of the other case is analogous. \qed

Let $(\FF^{(i)}, \GG^{(i)}, P^{(i)}, V^{(i)})$ be the quadruple obtained by iterating the above procedure $i$ times.

\bigskip

{\bf Summary of output of procedure.} The procedure applied to $(\FF^{(i)}, \GG^{(i)}, P^{(i)}, V^{(i)})$ with $V^{(i)}:\cong[n-i]$ either stops or yields
$(\FF^{(i+1)}, \GG^{(i+1)}, P^{(i+1)}, V^{(i+1)})$ with $V^{(i+1)}:\cong[n-i-1]$.  Moreover,

i) $p(\FF^{(i+1)})p(\GG^{(i+1)})> a \, p(\FF^{(i)}) p(\GG^{(i)})$ and $P^{(i+1)}\in \{P^{(i)}-1, P^{(i)}\}$ or

ii) $p(\FF^{(i+1)})p(\GG^{(i+1)})> \eps \, p(\FF^{(i)}) p(\GG^{(i)})$ and $P^{(i+1)}= (P^{(i)}-1) \cup P^{(i)}$.

Observe that ii) can occur at most $l$ times.  Indeed, if $P^{(i+1)} \supset (P^{(i)}-1) \cup P^{(i)}$ happens for $l$ values of $i$, then (since $l([n]\setminus P)=l$) the resulting set of forbidden intersection sizes consists of all nonnegative integers.  Consequently one of the resulting families would have to be empty.  This however means that all of (\ref{star}), (\ref{2star}), (\ref{3star}) fail and hence the algorithm stops.
 \qed

We will begin the procedure with $(\FF, \GG)=(\FF^{(0)}, \GG^{(0)})=(\AA, \BB)$ and distinguish two cases.

Suppose that the procedure does not stop until $i=n$. Then for each $i\le n$ we have
$$p(\FF^{(i+1)})p(\GG^{(i+1)})\ge a\, p(\FF_{i})p(\GG_{i}) \quad \hbox{ or } \quad p(\FF^{(i+1)})p(\GG^{(i+1)})\ge \eps\, p(\FF_{i})p(\GG_{i}).$$
Since ii) occurs at most $l$ times and  $\eps<a$ we obtain
\begin{equation} \label{i=n} 1\ge p(\FF_{n})p(\GG_{n})\ge a^{n-l}\eps^l p(\AA)p(\BB).\end{equation}
Consequently,
$$|\AA||\BB|=p(\AA)p(\BB)4^n<\frac{4^n}{a^{n-l}\eps^l}<
K\left(\frac{4}{a}\right)^n$$
where $K=(a/\eps)^l<(2/\eps)^l$.  As $\eps=2/10^4$, and $4/a<2.631$, this implies that
$$|\AA||\BB|<\left(\frac{4}{a}\right)^n  \left(\frac{2}{\eps}\right)^l  < 2.631^n \times 10^{4l}.$$

Suppose now that the procedure stops at some stage $i <n$. 
Then the Claim shows that  for all $v \in V^{(i)}:\cong[n-i]$ either
$$ p(\FF^{(i)}_1(v)) \le \sqrt{\eps} p(\FF^{(i)}) \quad \hbox{ or } \qquad
p(\GG^{(i)}_1(v)) \le \sqrt{\eps} p(\GG^{(i)}).$$
We will prove  a similar bound for $|\AA||\BB|$ in this case. 
Let
$$W_{\FF}=\{v \in V^{(i)}: p(\FF^{(i)}_1(v)) < \eps^{1/2}\, p(\FF^{(i)})\}
\quad \hbox{ and } \quad
W_{\GG}=V^{(i)}\setminus W_{\FF}.$$
Observe that $p(\GG^{(i)}_1(v)) < \sqrt{\eps} p(\GG^{(i)})$ for all $v
\in W_{\GG}$.

Let us now focus on $\FF^{(i)}$ and $W_{\FF}$. Consider the bipartite graph $H_{\FF}$ with bipartition $W_{\FF}$ and $\FF^{(i)}$ where $v \in W_{\FF}$ is joined to $A \in \FF^{(i)}$ if $A\setminus \{v\} \in \FF_1(v)$. Then the degree $\hbox{deg}_{H_{\FF}}(v)$ of  $v \in W_{\FF}$ in $H_{\FF}$ is 
$$|\FF^{(i)}_1(v)|=p(\FF^{(i)}_1(v))2^{|V^{(i)}|-1}<
{\eps}^{1/2}p(\FF^{(i)})2^{|V^{(i)}|-1}
={\eps}^{1/2}\frac{|\FF|}{2^{|V^{(i)}|}}2^{|V^{(i)}|-1}
=\frac{1}{2}\eps^{1/2}|\FF^{(i)}|.$$
Let $x:= |W_{\FF}|$ and $s:=\lceil\frac{1}{2}\e^{1/2}x\rceil $. Counting the edges of $H_{\FF}$ in two different ways yields
\begin{equation} \label{xs}\sum_{S \in \FF^{(i)}} |S \cap W_{\FF}|= \sum_{S \in \FF^{(i)}} \hbox{deg}_{H_{\FF}}(S) =\sum_{v \in W_{\FF}} \hbox{deg}_{H_{\FF}}(v)< \frac{1}{2}\eps^{1/2}x|\FF^{(i)}|\le s|\FF^{(i)}|.
\end{equation}
For integers $a \ge b\ge 0$, write ${a \choose \le b}$ to denote $\sum_{i=0}^b {a \choose i}$. 

{\bf Claim.} $|\FF^{(i)}| \le (s+1) {x \choose \le s} 2^{|V^{(i)}|-x}$.

{\bf Proof of Claim.} Let us suppose for contradiction that
$|\FF^{(i)}| > (s+1) {x \choose \le s} 2^{|V^{(i)}|-x}$.  Write $\FF^{(i)}=\FF' \cup \FF''$ where
$$\FF'=\{ S \in \FF^{(i)}: |S \cap W_{\FF}|\le s\} \quad \hbox{ and } \quad \FF''=\FF^{(i)}\setminus \FF'.$$
Then, trivially
$$|\FF'| \le {x \choose \le s} 2^{|V^{(i)}|-x}< \frac{1}{s+1}|\FF^{(i)}|.$$
Consequently,
$$\sum_{S \in \FF^{(i)}} |S \cap W_{\FF}|\ge \sum_{S \in \FF''} |S \cap W_{\FF}|>
\frac{s}{s+1} |\FF^{(i)}| (s+1)=s|\FF^{(i)}|.$$
This contradicts (\ref{xs}) and the claim is therefore proved. \qed

Similarly,
$$|\GG^{(i)}|\le (t+1) {y \choose \le t}2^{|V^{(i)}|-y}$$
where  $y:=|W_{\GG}|$ and $t:= \lceil\frac{1}{2}\e^{1/2}y\rceil $.
Note that
$x+y=|W_{\FF}|+|W_{\GG}|= |V^{(i)}|=n-i$. Putting this together we get
\begin{align}
|\FF^{(i)}||\GG^{(i)}| &\le (s+1)(t+1){x \choose \le s}{y \choose \le t}2^{|V^{(i)}|-x}2^{|V^{(i)}|-y} \n \\
&\le (s+1)^2(t+1)^2 {x \choose  s}{y \choose t}2^{|V^{(i)}|-x}2^{|V^{(i)}|-y}. \n \end{align}
The inequality $ab\le (a+b)^2/4$ for positive reals $a,b$ yields  
$$(s+1)^2(t+1)^2 \le \left(\frac{s+t+2}{2}\right)^4 \le \left(\frac14 
\e^{1/2}(n-i)+1\right)^4=\left(\frac{1}{16}\e(n-i)^2+\frac12\e^{1/2}(n-i)+1\right)^2.$$
Since $\e=2/10^4$ and $n-i \ge 1$, 
$$1 < \frac{10^5}{3} \e (n-i)^2 \quad \hbox{ and } \quad
\frac12 \eps^{1/2} (n-i) < \frac{10^5}{3} \e (n-i)^2.$$
Consequently, $(s+1)^2(t+1)^2 < 10^{10} \e^2 (n-i)^4$, and we therefore have 
\begin{align}
|\FF^{(i)}||\GG^{(i)}|&\le 10^{10}\eps^2(n-i)^4 {x \choose  s}{y \choose t}2^{|V^{(i)}|-x}2^{|V^{(i)}|-y} \n \\
&\le 10^{10}\eps^2(n-i)^4 \left(\frac{xe}{s}\right)^s \left(\frac{ye}{t}\right)^t 2^{|V^{(i)}|-x}2^{|V^{(i)}|-y} \n \\
&\le
10^{10}\eps^2(n-i)^4\left(\frac{2e}{\eps^{1/2}}\right)^{\frac{1}{2}\eps^{1/2}(x+y)}2^{2n-2i-x-y}\n \\
&=10^{10}({\eps^{1/2}}(n-i))^4\left(\frac{2e}{\eps^{1/2}}\right)^{\frac{1}{2}\eps^{1/2}(n-i)}
2^{n-i}. \n\end{align}
Set
$$\sigma=\frac{4\log_2( \eps^{1/2}(n-i))}{n-i}+
\frac{1}{2}\eps^{1/2} \log_2(2e/\eps^{1/2})$$ and suppose that the algorithm stops at stage $i$ with $(\FF^{(i)}, \GG^{(i)})$.
Then we have shown above that $|\FF^{(i)}||\GG^{(i)}|<10^{10}2^{(n-i)(1+\sigma)}$. As in (\ref{i=n}), we obtain
$$p(\AA)p(\BB)a^{i-l}
\eps^l<p(\FF^{(i)})p(\GG^{(i)})<
\frac{10^{10}2^{(n-i)(1+\sigma)}}{2^{2(n-i)}}=10^{10} 2^{(n-i)(\sigma-1)}.$$
Consequently,
\begin{equation} \label{11}|\AA||\BB| \le 4^n p(\AA)p(\BB)\le
4^n \frac{10^{10} 2^{(n-i)(\sigma-1)}}{a^{i-l}\eps^l}=\frac{10^{10} 4^n}{2^{n(1-\sigma)}}
\left(\frac{2^{(1-\sigma)}}{a}\right)^i\left(\frac{a}{\eps}\right)^l.\end{equation}

We now provide an upper bound for $\sigma$.  Consider the real valued function $f(m)= 4\log_2(\eps^{1/2}m)/m$ where $1\le m \le n$. Easy calculus shows that $f(m)$ is maximized when $m=e\eps^{-1/2}$, where its value is $(4/e \ln 2)\eps^{1/2}$.
Using $\eps=2/10^4$ and $a<6-2\sqrt{5}<1.53$, we obtain
$$\sigma \le \frac{4 \eps^{1/2}}{e\ln 2} + \frac{1}{2} \eps^{1/2} \log_2(2e/\eps^{1/2})<0.1<0.38<1-\log_2 a.$$
Consequently, $2^{1-\sigma}>a$ and
 the RHS in (\ref{11}) is upper bounded by
$$ 10^{10} 4^n \frac{1}{2^{(1-\sigma)(n-i)} a^i} \left(\frac{a}{\eps}\right)^l< 10^{10}
\left(\frac{4}{a}\right)^n \left(\frac{a}{\eps}\right)^l.$$ This implies that
$$|\AA||\BB|< 10^{10}\left(\frac{4}{a}\right)^n\left(\frac{2}{\eps}\right)^l < 2.631^n \times 10^{4l+10}.$$

{\bf Remark.}  As mentioned earlier, we have not optimized the value of $\eps$ in our proof, indeed the inequality $0.1<0.38$ above reflects the slack in our calculations.

\subsection{The second bound}
In this section we prove that $|\AA||\BB|<2^{n+2l\log^2 n}$.  We will use a very general result of Sgall~\cite{Sg} and show that it can be applied to our setting.  To describe the result of Sgall, we need some definitions.

\begin{definition} \label{1hdef}  Say that a function $h : \mathbb{N}^{<\infty} \rightarrow \mathbb{N}\cup \{\infty\}$ is a {\em height function} if the following four properties hold:

(A1) $h(L)=0$ if and only if $L=\emptyset$,

(A2) if $h(L)<\infty$ and $L' \subset L$, then $h(L')\le h(L)$,

(A3) if $h(L)<\infty$ and $L' \subset L-1$, then $h(L')\le h(L)$,

(A4) if $h(L), h(L')\le s <\infty$, then either  $h(L' \cap L)\le s-1$ or
$h(L' \cap (L-1))\le s-1$.
\end{definition}

\begin{definition} \label{sig}
Given a family $\AA$ and a set $B$, define the {\em signature} of $B$ to be the set
$$L_{B}^{\AA}=\{|A \cap B|: A \in \AA\}.$$
\end{definition}

\begin{definition} \label{hf}
A pair of families $(\AA, \BB)$ has height  $s$ if there is a height function $h$ such that for all $B \in \BB$ we have $h(L_{B}^{\AA})\le s$. 
\end{definition}

We now state Sgall's theorem.

\begin{theorem} \label{sgall} {\bf (Sgall~\cite{Sg})}
Suppose that $(\AA, \BB)$ is a pair of families on $V$, $|V|=n$ and $(\AA, \BB)$ has height  $s \le n+1$. Then
$$|\AA||\BB|\le 2^{n+s-1}{n \choose s-1}.$$
\end{theorem}

If $(\AA, \BB)$ is $M$-intersecting, and there is a height function $h$ with $h(M) \le s$, then  $(\AA, \BB)$ has height  $s$.  Indeed, this holds because $M =\cup_{B \in \BB} L_B^{\AA}$ and by (A2) we have $h(L_B^{\AA}) \le h(M)$ for every $B \in \BB$.
 
In order to prove the second bound in Theorem \ref{main} it will therefore be sufficient to  define a height function $h$ such that $h(M) \le 1+2l(M)\log n$ as long as $l(M)<n/(2\log n)$ (this will ensure that $h(M) \le n+1$).  Then, if $l(M)=l$, and $n$ is sufficiently large, Theorem \ref{sgall} immediately gives
$$|\AA||\BB|\le 2^{n+2l\log n}{n \choose 2l\log n} < 2^{n+2l\log^2 n}.$$

{\bf Definition of height function $h$.} The height function $h$ is defined recursively. First, let $h(\emptyset)=0$.  Now suppose that $L\ne \emptyset$ and $h$ has been defined on all sets  with size less than $|L|$.  Then
\begin{equation}\label{hdef}
h(L)=1+\max\{h(L \cap (L+1)), \max_{M \in T(L)} \,\min\{h(L \cap M),\, h(L \cap (M-1))\}\},\end{equation}
where
$$T(L)=\{M\subset \mathbb{N}^{<\infty}: M\not\in \{L, L+1\} \hbox{ and } 0<|M|\le |L|\}.$$

  An easy consequence of the above definition is that $h(L)=1$  if $|L|=1$.

Another easy consequence of this definition is that $h$ satisfies (A2).  Indeed,
let us prove (A2) by induction on $|L|$.  The result clearly holds if $|L|=1$ so
let $|L| \ge 2$ and $L' \subsetneq L$.  First suppose that $h(L')=1+h(L' \cap (L'+1))$. Since
$|L \cap (L+1)|<|L|$ and $L' \cap (L'+1) \subset L \cap (L+1)$, we can apply induction to get
$$h(L') =1+h(L' \cap (L+1))\le 1+ h(L \cap (L+1)) \le h(L)$$
where the last inequality holds by the definition of  $h(L)$.

We may now suppose that there is an $M \in T(L')\subset T(L)$ with
$$h(L')=1+\min\{h(L' \cap M), \, h(L' \cap (M-1)) \}.$$
Since $|M|\le |L'|<|L|$, we have $|L \cap M|<|L|$ and $|L \cap (M-1)|<|L|$. As $L' \cap M \subset L\cap M$ and $L' \cap (M-1) \subset L \cap (M-1)$ we have
by induction
$$1+\min\{h(L' \cap M), \, h(L' \cap (M-1)) \}\le 1+\min\{h(L \cap M), \, h(L \cap (M-1)) \}.$$
Since $M \in T(L') \subset T(L)$, the RHS above is at most
$$ 1+\max_{M \in T(L)}\min\{h(L \cap M), \, h(L \cap (M-1)) \}\le h(L)$$
and therefore $h(L') \le h(L)$.

Having shown that $h$ satisfies (A2) allows us to give a slightly simpler expression for $h$ as follows.
 
 \begin{proposition} The height function $h$ satisfies  $h(\emptyset)=0$, and if $|L|>0$, then
\begin{equation} \label{newh}h(L)=1+ \max_{M \in S(L)} \,\min\{h(L \cap M),\, h(L \cap (M-1))\},\end{equation}
where
$$S(L)=\{M\subset \mathbb{N}^{<\infty}: M\ne L \hbox{ and } 0<|M|\le |L|\}.$$
\end{proposition}

\proof Let $t=\max_{M \in T(L)} \,\min\{h(L \cap M),\, h(L \cap (M-1))\}$ so that by definition, $$h(L)=1+\max\{h(L \cap (L+1)), t\}.$$   By A2, $h(L \cap (L+1)) \le h(L)$, so
$$\min\{ h(L \cap (L+1)), h(L)\} = h(L \cap (L+1)).$$
Consequently, 
\begin{equation} \label{hl}
h(L)=1+\max\{ \min\{ h(L \cap (L+1)), h(L)\}, t\}\end{equation}
Since $S(L)=T(L) \cup \{L+1\}$ the RHS of (\ref{newh}) equals the RHS of (\ref{hl}).
\qed

To finish the proof, it suffices to prove the following two propositions.

\begin{proposition} \label{p1}
The function $h$ defined above is a height function
\end{proposition}

\begin{proposition} \label{p2}
If $M \subset [n]$ and $l(M)= l$, then $h(M) \le 1+2l\log n$.
\end{proposition}

We now prove each of these propositions.

{\bf Proof of Proposition \ref{p1}.} We must show that (A1)--(A4) hold.  Clearly (A1) holds by definition and we have already shown that (A2) holds.  Let us now prove that (A3) and (A4) hold.  In what follows we will use the expression for $h(L)$ given in (\ref{newh}).

{\bf(A3)} We will prove (A3) by induction on $|L|$.  The result clearly holds if $|L|=1$ so
let $|L| \ge 2$. By (\ref{newh}), there is an $M \in S(L-1)$ with
$$h(L-1)=1+\min\{h((L-1) \cap M),\, h((L-1) \cap (M-1))\}.$$
Since $M \ne L-1$, we have $|L \cap (M+1)| < |L|$.  W will distinguish two cases.

a) If $M\ne L$ we  have $|L \cap M|<|L|$.
Hence by induction $$h((L-1) \cap M) \le h(L\cap (M+1)).$$
Similarly, we also have
$$h((L-1) \cap (M-1)) \le h(L \cap M).$$  Since $M \in S(L-1)$, we have $M+1 \in S(L)$.  Thus
$$h(L-1)\le 1+\min\{h(L \cap (M+1)),\, h(L \cap M)\}\le h(L)$$
and hence A3 holds on the assumption $M\ne L$.

b) If $M=L$, then by (A2), $h((L-1) \cap M)\le h((L-1) \cap (M-1))$, and hence
$$h(L-1)=1+h((L-1) \cap L).$$
Furthermore, $|L\cap (L+1)| < |L|$  so by induction, $h((L-1) \cap L)\le h(L \cap (L+1))$. Consequently,
\begin{align}
h(L-1)&=1+\min\{h((L-1) \cap M), h((L-1) \cap (M-1))\} \notag \\
&\le  1+\min\{h(L \cap (M+1)), h(L \cap M)\}\le h(L) \notag .
\end{align}

\medskip

{\bf (A4)} Let $L, L'$ be given with $h(L), h(L')\le s$.  First let us consider the case that $L=L'$.
In that case, let us assume for contradiction that $h(L)\ge s$ and $h(L \cap (L-1)) \ge s$. Now let $M=L+1 \in S(L)$ and thus $h(L \cap (M-1))=h(L)\ge s$.  Then by (A3), $h(L \cap M) \ge h(L \cap (L-1))\ge s$.  Consequently, we have the contradiction
$$h(L) \ge 1+\min\{h(L \cap M), h(L \cap (M-1))\} \ge s+1.$$
So we henceforth assume that $L \ne L'$.
Now by (\ref{newh})
\begin{equation} \label{e1}
\max_{M \in S(L)} \min \{ h(L \cap M), h(L \cap (M-1))\}\le s-1
\end{equation}
 and
 \begin{equation} \label{e2}
\max_{M'\in S(L')} \min \{ h(L' \cap M'), h(L' \cap (M'-1))\}\le s-1
\end{equation}

Let use first suppose that $|L| \le |L'|$ and put $M'=L$  in (\ref{e2}).  Notice that $L \ne L'$ and $|L|\le |L'|$ yield $M' \in S(L')$.  Then (\ref{e2}) gives precisely what we want:
$$\min \{h(L' \cap L), \, h(L' \cap (L-1))\} \le s-1.$$
Next suppose that $|L'|<|L|$ and put $M=L'+1$ in (\ref{e1}); since $|L'|<|L|$ we have $L'+1 \ne L$ and hence $M \in S(L)$. Then (A3) implies that
$$h(L' \cap (L-1)) \le h((L'+1) \cap L)=h(L \cap M).$$
Finally, (\ref{e1}) yields
$$\min\{h(L' \cap L), h(L' \cap (L-1))\}\le \min\{h(L\cap (M-1)), h(L \cap M)\} \le s-1.$$

This completes the proof of the proposition. \qed
\bigskip

{\bf Proof of Proposition \ref{p2}.} Let $W \subset [n]$ and $l(W)=l$.  We are going to show that $h(W) \le 1+2l\log n$.
We will prove by induction on $|W|$ that if $l(W)=l$, then 

\begin{equation} \label{starr}
\left(\frac{2l}{2l-1}\right)^{h(W)-1} \le |W|.\end{equation}
 The result is trivial if $|W|=l$ ($|W|<l$ is impossible), since in this case (\ref{newh}) implies that $h(W)\le l$. The identity $(2l/(2l-1))^{l-1}\le l$ is now easily checked. So assume that $|W|>l$. For the induction step, we have
\begin{equation} \label{Vstar}
h(W) =1+ \min \{h(W \cap M), h(W \cap (M-1))\}\end{equation}
for some $M \in S(W)$.
We may assume that $W$ is critical, namely that if $W' \subsetneq W$, then  $h(W')<h(W)$.  This is because if $W$ is not critical, then there is some critical $W'\subsetneq W$ with
$h(W')=h(W)$, and if we have proved the result for critical sets, then
$$|W|\ge |W'|\ge (2l/(2l-1))^{h(W')-1}=(2l/(2l-1))^{h(W)-1}.$$
Consider the two sets $W \cap M$ and $W \cap (M-1)$.
To every element $x \in W \cap (M-1)$ associate the element $x+1 \in M$ (this is clearly an injection). Since $l(W)=l$, we can write $W$ as a union of disjoint intervals each of length at most $l$.  The first element of each of these intervals cannot belong to $W \cap (W+1)$, and there are at least $|W|/l$ such elements.  Consequently, $|W \cap (W+1)|\le (1-1/l)|W|$, and so
\begin{align}
|W \cap M| + |W \cap (M-1)|&= |W \cap M| + |(W +1) \cap M| \notag\\
&\le |W \cap (W +1)|+|M| \notag \\
&\le |M| + (1-1/l)|W| \le (2-1/l)|W|.\notag \end{align}
So either $|W \cap M|<|W|(2l-1)/(2l)$ or $|W \cap (M-1)|\le |W|(2l-1)/(2l)$.
 Suppose the former holds. Since $M\ne W$ and $W$ is critical, $h(W \cap M)\le h(W)-1$.
  On the other hand by (\ref{Vstar}) $h(W)-1\le h(W \cap M)$ and thus $h(W)-1=h(W \cap M)$.
     So by induction on $|W \cap M|$,
\begin{equation} \label{V2star}(2l/(2l-1))^{h(W)-1}=(2l/(2l-1))^{h(W \cap M)} \le (2l/(2l-1))|W \cap M|\le |W|.\end{equation}
 Next suppose that $|W \cap (M-1)| \le |W|(2l-1)/(2l)$.  In this case $W \ne M-1$ must hold and since $W$ is critical, $h(W \cap (M-1))=h(W)-1$. Now we apply induction  as in (\ref{V2star}) with $W \cap M$ replaced by $W \cap (M-1)$. In either case we obtain
 $|W| \ge (1+1/(2l-1))^{h(W)-1}$ or equivalently,
$$h(W) \le 1+\frac{\log|W|}{\log\left(1+\frac{1}{2l-1}\right)}.$$
as long as $l(W)=l$. This completes the induction proof of (\ref{starr}).

Since $\log(1+x)>x-x^2/2$ for  $|x|<1$, we have
$$\log\left(1+\frac{1}{2l-1}\right)>\frac{1}{2l-1}-\frac{1}{2(2l-1)^2}
=\frac{4l-3}{2(2l-1)^2}\ge \frac{1}{2l}.$$
Inserting this above yields
$h(W) < 1+2l\log n$ as required.
\qed

\section{Proof of Theorem \ref{misint}}
In this section we present the short proof of Theorem \ref{misint}.

Our main tool is the following result which follows from the Frankl-Wilson Theorem~\cite{FW}.

\begin{theorem} \label{misint3} Let
$n>k>2t$. Suppose that $\BB$ is a $\{t\}$-omitting family of $k$-element subsets of
$[n]$. If
$$
gcd\left({k-1 \choose k-t-1}, {k-2 \choose k-t-1}, \ldots, {k-t
\choose k-t-1}\right)>1, \qquad (*)
$$
then $|\BB| \le {n \choose k-t-1}.$
\end{theorem}
It is easy to see that ($*$) holds if, for example, $k-t$ is prime.

We will also use the result of Baker-Harman-Pintz \cite{BHP} which states that for all $s$ sufficiently large, there is a prime in the interval $(s-s^{0.525}, s)$.

{\bf Proof of Theorem \ref{misint}.}
We will omit floor and ceiling symbols. Let $\gamma=0.525$, $0<\eps<1/5$  and assume that a $(t, t+n^{\gamma})$-omitting family $\AA \subset 2^V$ is given, and $n>n_0(\eps)$.
Write $\AA_k$ for the family of those
subsets of $\AA$ of size exactly $k$. Then

$$|\AA|=\sum_{k \le 2t}|\AA_k| + \sum_{2t < k\le (n+t)/2}|\AA_k|
+\sum_{k>(n+t)/2}|\AA_k|.$$

Each term in the first  summation is bounded by ${n \choose 2t}<{n
\choose (n+t)/2}$ since $t<n/5$. Each term in the last summation is
clearly bounded by ${n \choose (n+t)/2}$.

Now consider $\AA_k$ with $2t<k \le (n+t)/2$. Since $\AA$ is $\{t\}$-omitting,
$\AA_k$ is $\{t'\}$-omitting for every $t' \in (t, t+n^{\gamma})$.  By the result of \cite{BHP}, since $k-t>t>\eps n$ and $n>n_0(\eps)$ we can find a prime $p \in  (k-t-(k-t)^{\gamma} , k-t)$. So $p=k-t_k$ where $t \le t_k \le t+n^{\gamma}$ and $\AA$ is $\{t_k\}$-omitting.
Now apply Theorem \ref{misint3} to bound each term in the second
summation by ${n \choose k-t_k-1}$. Since $t_k\ge t$, the bounds in the second summation
are at most
$${n \choose t}, \ldots, {n \choose (n-t)/2-1}.$$ As $t<n/5$, we have
$t<(n-t)/2<n/2$, and each of these terms is less than ${n \choose
(n-t)/2}$.  Thus we get
$$|\AA|\le n {n \choose (n-t)/2}=n {n \choose (n+t)/2}$$
and the proof is complete. \qed

\section{Acknowledgments}
We thank Sundar Vishwanathan for reading an earlier version of this manuscript and suggesting some simplifications in the presentation. Thanks also to Zoltan F\"uredi for informing us about \cite{FF}, Mathias Schacht for informing us about \cite{BCW} and \cite{KG}, and a referee for helpful remarks about the presentation.

\end{document}